\documentclass [12pt,a4paper] {article}
\usepackage[affil-it]{authblk}

\usepackage [T1] {fontenc}
\usepackage {amssymb}
\usepackage {amsmath}
\usepackage {amsthm}
\usepackage{hyperref}
\usepackage{tikz}

\usetikzlibrary{arrows,decorations.pathmorphing,backgrounds,positioning,fit,petri}

\usetikzlibrary{matrix}

\DeclareMathOperator {\td} {td}
\DeclareMathOperator {\U} {U}
\DeclareMathOperator {\Th} {Th}
\DeclareMathOperator {\ldim} {l.dim}
\DeclareMathOperator {\tp} {tp}
\DeclareMathOperator {\qftp} {qftp}
\DeclareMathOperator {\ord} {ord}
\DeclareMathOperator {\acl} {acl}
\DeclareMathOperator {\dcl} {dcl}

\DeclareMathOperator {\dif} {dif}

\DeclareMathOperator {\DCF} {DCF}
\DeclareMathOperator {\ACF} {ACF}
\DeclareMathOperator {\DR} {DR}
\DeclareMathOperator {\DF} {DF}
\DeclareMathOperator {\MR} {MR}
\DeclareMathOperator {\MD} {MD}
\DeclareMathOperator {\D} {D}

\DeclareMathOperator {\alg} {alg}

\DeclareMathOperator {\E} {E}

\DeclareMathOperator {\EDiag} {EDiag}
\DeclareMathOperator {\Aut} {Aut}

\theoremstyle {definition}
\newtheorem {definition}{Definition} [section]
\newtheorem{example} [definition] {Example}
\newtheorem* {claim} {Claim}

\theoremstyle {plain}
\newtheorem {question} [definition] {Question}
\newtheorem {lemma} [definition] {Lemma}
\newtheorem {theorem} [definition] {Theorem}
\newtheorem {proposition} [definition] {Proposition}
\newtheorem {corollary} [definition] {Corollary}

\theoremstyle {remark}
\newtheorem {remark} [definition] {Remark}

\makeatletter
\newcommand {\forksym} {\raise0.2ex\hbox{\ooalign{\hidewidth$\vert$\hidewidth\cr\raise-0.9ex\hbox{$\smile$}}}}
\def\@forksym@#1#2{\mathrel{\mathop{\forksym}\displaylimits_{#2}}}
\def\forkind{\@ifnextchar_{\@forksym@}{\forksym}}
\newcommand {\nforkind} {\not \forkind}
\makeatother

\makeatletter
\def\@maketitle{%
  \newpage
  \null
  \vskip 2em%
  \begin{center}%
  \let \footnote \thanks
    {\Large\bfseries \@title \par}%
    \vskip 1.5em%
    {\normalsize
      \lineskip .5em%
      \begin{tabular}[t]{c}%
        \@author
      \end{tabular}\par}%
    \vskip 1em%
    {\normalsize \@date}%
  \end{center}%
  \par
  \vskip 1.5em}
\makeatother

\begin {document}

\title{Definability of Derivations in the Reducts of Differentially Closed Fields}

\author{Vahagn Aslanyan\footnote{E-mail: \texttt{vahagn.aslanyan@gmail.com}}}
\affil{Mathematical Institute, University of Oxford, Oxford, OX2 6GG, UK}

\maketitle

\begin {abstract}
Let $\mathcal{F}=(F;+,\cdot,0,1,\D)$ be a differentially closed field. We consider the question of definability of the derivation $\D$ in reducts of $\mathcal{F}$ of the form $\mathcal{F}_{R}=(F;+,\cdot,0,1,P)_{P \in R}$ where $R$ is some collection of definable sets in $\mathcal{F}$. We give examples and non-examples and establish some criteria for definability of $\D$. Finally, using the tools developed in the paper we prove that under the assumption of inductiveness of $\Th(\mathcal{F}_{R})$ model completeness is a necessary condition for definability of $\D$. This can be seen as part of a broader project where one is interested in finding Ax-Schanuel type inequalities (or predimension inequalities) for differential equations.

\end {abstract}

\begin{small}

\noindent \textbf{2010 MSC:} 03C10, 03C60, 12H05, 12H20, 13N15.

\noindent \textbf{Keywords:} Model theoretic algebra, differentially closed field, reduct, abstract differential equation, definable derivation.
\end{small}

\section{Introduction}

For a differentially closed field $\mathcal{F}=(F;+,\cdot,0,1,\D)$ we consider its reducts of the form $\mathcal{F}_{R}=(F;+,\cdot,0,1,P)_{P \in R}$ where $R$ is some collection of definable sets in $\mathcal{F}$. Our main problem is to understand when the derivation $\D$ is definable in $\mathcal{F}_{R}$. Ideally, we would like to find a dividing line for definability of $\D$ like local modularity in the problem of recovering the field structure in the reducts of algebraically closed fields (see the discussion below).

\begin{question}\label{Question-definability}
When is $\D$ definable in the reduct $\mathcal{F}_{R}$?
\end{question}

As we will see when $\D$ is definable it is definable with using just one parameter, namely an element $t \in F$ with $\D t =1$. So it is more convenient to add $t$ to our language as a constant symbol and work in the reducts of $\mathcal{F}=(F;+,\cdot,0,1,t,\D)$ (we do this starting from Section \ref{some examples}). So, we will assume for simplicity that the sets from $R$ are $0$-definable in this language and also we will be interested in $0$-definability of $\D$.

Note that one could be tempted to ask a more general
question of whether there is a derivation definable in the reduct. But in that case such a derivation will also be definable in the differentially closed field $\mathcal{F}$. Since it is known that any such derivation is of the form $a \cdot \D$ for some $a\in F$, it is no loss of generality if we restrict our attention to definability of $\D$ only. Another point is that we can assume that $R$ is finite since any possible definition of $\D$ can contain only finitely many relations from $R$.

This is by nature a classification problem. We do not have a comprehensive solution yet, but we give some partial answers to our question, and draw some conclusions based on our analysis. We will not pose any explicit conjectures, but one may nevertheless expect intuitively that definability of $\D$ is very rare, i.e. in most cases it is not definable. In other words, our general expectation is that for ``generic'' reducts $\D$ is not definable.

The motivation to consider this kind of problem comes from two independent sources. Firstly, the analogous problem for pure fields, that is, recovering the field structure from reducts of algebraically closed fields or from non-locally modular strongly minimal sets in general, is a well studied question in model theory of fields and Zariski geometries. It was initiated by Zilber's famous ``Trichotomy conjecture'' and is still not entirely resolved. It has been (and still is) a topic of active research during the past few decades and proved to be very useful. Zariski geometries, introduced by B. Zilber and E. Hrushovski, are structures where that theory works ideally. For more details on this we refer the reader to \cite{Zariski,Rabin,Hasson,Mar-dif}.

Secondly, this problem turns out to be related to the existence of an ``Ax-Schanuel type theorem'' for a given differential equation $E(x,y)$ (in this case we will work in the reduct $\mathcal{F}_{E}=(F;+,\cdot, 0,1,E)$ with $R=\{ E\}$). Let us briefly explain what we mean by this.

James Ax has proved the following analogue of Schanuel's conjecture in differential setting (\cite{Ax}). Let $K$ be a differential field and $C$ be its field of constants. Let also $(x_1,y_1),\ldots,(x_n,y_n)$ be non-constant solutions to the exponential differential equation $\D x=\frac {\D y} {y}$ in $K$. Then
$$\delta(x_1,\ldots,x_n):= \td_CC(x_1,\ldots,x_n,y_1,\ldots,y_n)-\ldim_{\mathbb{Q}}(x_1,\ldots,x_n/C) \geq 1,$$
where $\td$ stands for the transcendence degree and $\ldim$ stands for the linear dimension (modulo $C$) as a vector space. This inequality is now known as the Ax-Schanuel inequality.
The function $\delta$ here is a \textit{predimension} function in the sense of Hrushovski (\cite{Hru}). Thus the Ax-Schanuel inequality is a \textit{predimension inequality}. This property gives a good understanding of the exponential differential equation. In particular one can consider the corresponding reduct. Then the first order theory of the reduct is axiomatised by axioms of algebraically closed fields, functional equation(s), an axiom scheme for the Ax-Schanuel inequality and the strong existential closedness axiom scheme (see \cite{Zilb1, Kirby-semiab}). This is exactly the axiomatisation that one obtains after carrying out a Hrushovski construction with the above predimension function. Thus, the reduct here is reconstructed by a Hrushovski construction. Zilber calls such predimension inequalities \emph{adequate}.

After realising this one can ask whether it is possible to do something similar for other differential equations. One therefore poses a problem whether for a given differential equation there is an ``Ax-Schanuel type'' inequality (or a predimension inequality). It is useful to classify differential equations with respect to this property, i.e. whether there is an ``adequate'' predimension inequality or not. If there is one, then one will know the complete theory of the equation. One of the recent developments in this direction is the establishment of an Ax-Schanuel type inequality for the $j$-function by Jonathan Pila and Jacob Tsimerman (\cite{PilaTsim}). For details on Schanuel's conjecture and the Ax-Schanuel inequality (and its generalised versions) see \cite{Zilb1,Zilb2,Zilb3,Zilb4,Kirby-semiab}. For Hrushovski constructions and predimensions we refer the reader to \cite{Hru,Wagner}.

We are not going to consider these questions in this paper, but let us see how this problem is related to definability of $\D$ in the corresponding reduct. The idea is that definability of a derivation would imply that there is no ``non-trivial'' adequate predimension inequality for the given differential equation. Roughly speaking, if $\D$ is definable then the problem is reduced to finding an Ax-Schanuel type inequality for the equation $y=\D x$. But one can argue that there is no such non-trivial inequality for the latter equation. 
We will support this viewpoint by a result in the last section. Indeed, as we will see if $\D$ is definable and the theory of the reduct is inductive then it must actually be model complete. But  Hrushovski constructions yield inductive theories in nice examples, so it is a reasonable condition. Thus, assuming the reduct is inductive, definability of $\D$ implies model completeness which can be used to show that a possible adequate predimension must be trivial in some sense. Of course, these statements are pretty vague and we presented them here just to give a basic idea about the connection of those two questions. 
We limit ourselves to these explanations and refer the reader to \cite{Aslanyan-thesis} for more details. 

Let us briefly outline the paper. After giving the necessary preliminaries in Section \ref{Preliminaries} we show in Section \ref{DefinableDerivations} that definable derivations in models of $\DCF_0$ are the trivial ones. Then we study the reducts of differentially closed fields from a general model theoretic point of view and establish some of their properties in Section \ref{ModelTheoreticProperties}. 
In Section \ref{some examples} we will see that if $E$ is a differential curve containing the graph of $\D$ then $\D$ is quantifier-free definable in $\mathcal{F}_{E}$. 

Furthermore, we will show in Section \ref{generic points} that the behaviour of $\D$ at generic points is enough to understand whether it is definable. Indeed we will prove that if for a generic element $a$ the Morley rank (in the reduct) of $\D a$ over $a$ is finite then $\D$ is definable (Theorem \ref{finiterank}). 
Using the results on generic points we will give further examples of differential equations that define $\D$ (Section \ref{Further examples}). Theorem \ref{equivalent-conditions} will sum up most of our results obtained up to that point giving a list of conditions equivalent to definability of $\D$ in the reducts.

The last section will be devoted to the question of model completeness of reducts that define $\D$. Namely, we will prove that if $\D$ is definable in $\mathcal{F}_{R}$ and $\Th(\mathcal{F}_{R})$ is inductive then this theory must in fact be model complete (Theorem \ref{inductive}). This will immediately imply that one cannot define $\D$ from the exponential differential equation $\D y = y\D x$.\\




\textbf{Acknowledgements.}
I would like to thank my supervisors Boris Zilber and Jonathan Pila for their support and numerous useful discussions. They have had great influence on this work and, in particular,  the idea that an adequate predimension inequality and definability of $\D$ together must imply model completeness is a result of many discussions with them. Therefore the formulation of Theorem \ref{inductive} (as a conjecture initially) is due to three of us.

I am also grateful to Ehud Hrushovski and Jonathan Kirby for reading this paper as part of my PhD thesis and making valuable comments.

Finally, I thank the referee for numerous useful remarks.

This research was supported by the University of Oxford Dulverton Scholarship.





\section{Preliminaries}\label{Preliminaries}

In this section we present basic definitions and facts about differential fields. For more details and proofs of the results stated here we refer the reader to
\cite{Mar-dif,Kap,Pillay1,Pillay2}.

We assume all rings that we deal with are commutative rings with identity and have characteristic zero.

The language of differential rings is $\mathfrak{L}_{\D}=\{+,\cdot, 0, 1, \D \}$. In this language we can axiomatise the theory of differential (rings) fields with the axioms of (rings) fields with two extra axioms stating that $\D$ is additive and satisfies Leibniz's rule, i.e. $\forall x,y~ \D(x+y)=\D x+\D y$ and $\forall x,y ~ \D(xy)=x\D y+y\D x$. The theory of differential fields of characteristic zero is denoted by $\DF_0$.

The \textit{field of constants} of a differential field $(F;+,\cdot,0,1,\D)$ is defined as the kernel of the derivation, i.e. $C_F=\{x\in F : \D x=0\}$. This is always a relatively algebraically closed subfield of $F$.

If $F$ is a differential field then the ring of \textit{differential polynomials} over $F$ is a differential ring extension defined as $F\{X\}=F[X,\D(X),\D^2(X),\ldots]$ with $\D(\D^n(X))=\D^{n+1}(X)$. Thus, differential polynomials are of the form $p(X,\D X,\ldots,\D^nX)$ where $p(X_0,\ldots,X_n)\in F[X_0,\ldots,X_n]$ is an algebraic polynomial over $F$. A \textit{differential rational function} over $F$ is the quotient of two differential polynomials over $F$. The field of all differential rational functions of $X$ over $F$ will be denoted by $F\langle X \rangle$. We can also consider differential polynomials in several variables, which are defined analogously. If $f(X_1,X_2,\ldots,X_n)$ is such a polynomial, then the equation $f=0$ is a \textit{differential equation} over $F$.

Further, for $F$ a differential field and $A\subseteq F$ a subset we denote by $\langle A \rangle$ or $\mathbb{Q}\langle A \rangle$ the differential subfield generated by $A$. If $K\subseteq F$ are differential fields and $A \subseteq F$ then $K\langle A\rangle$ is the differential subfield generated by $K$ and $A$. The algebraic subfield generated by $K$ and $A$ is denoted by $K(A)$. One can easily verify that $K\langle A\rangle=K(\{\D^n a: a\in A,~ n \in \mathbb{N} \})$.

The \textit{order} of $f$, denoted $\ord(f)$, is the biggest $n$ for which $\D^n(X)$ occurs in $f$. In this case the highest power of $\D^n(X)$ in $f$ is the \textit{degree} of $f$, written $\deg(f)$. In the case of polynomials of several variables we will write $\ord_{X_i}(f)$ for the order of $f$ with respect to $X_i$.

The theory $\DF_0$ has a model completion. It is called the theory of \textit{differentially closed fields} of characteristic zero. To axiomatise this theory we add the existential closedness axiom scheme: a differential field $(F;+,\cdot,0,1,D)$ is \textit{differentially closed} if for any non-constant differential polynomials $f(X)$ and $g(X)$ over $F$ with $\ord(g)<\ord(f)$ there exists $x\in F$ such that $f(x)=0$ and $g(x)\neq 0$. We let $\DCF_0$ denote the theory of differentially closed fields of characteristic $0$. It immediately follows from the definition that differentially closed fields are algebraically closed (in the field theoretic sense). Hence, the field of constants is algebraically closed as well.

Suppose $K\subseteq F$ are two models of $\DF_0$. For an element $a\in F$ one defines the \textit{differential rank} (or \textit{dimension} or \textit{order}) of $a$ over $K$, denoted $\DR(a/K)$ (or $\dim(a/K)$ or $\ord(a/K)$), as the transcendence degree of $K\langle a\rangle$ over $K$. If it is finite, say $n$, then there is a differential polynomial $f(X)\in K\{X\}$ of order $n$ with $f(a)=0$. If $f$ is the \textit{simplest} among such polynomials, i.e. the pair $(\ord(f),\deg(f))$ is minimal with respect to the lexicographical order, then it is called \textit{the minimal polynomial} of $a$ over $K$. This polynomial must be irreducible. The elements $a,\D a,\ldots,\D^{n-1}a$ are algebraically independent, while $a,\D a,\ldots,\D^na$ are algebraically dependent over $K$. In this case $a$ is called \textit{differentially algebraic} over $K$, otherwise it is called \textit{differentially transcendental} over $K$. In the latter case $\DR(a/K)$ is defined to be $\omega$.

Suppose $K\models \DF_0$ and $K\subseteq F$ is a differentially closed extension of $K$. Then for any element $a\in F$ the following inequality holds
$$\U(a/K) \leq \MR(a/K)\leq \DR(a/K),$$
where $\U(a/K)$ stands for the $\U$-rank and $\MR(a/K)$ stands for the Morley rank of $a$ over $K$. Moreover, $a$ is differentially transcendental over $K$ if and only if $\U(a/K)=\MR(a/K)=\DR(a/K)=\omega$. In this case $a$ is called \textit{generic} over $K$ (if we omit $K$ then it means $a$ is generic over the empty set or, equivalently, over the prime differential subfield).

There is a unique complete type of a differentially transcendental element (over a subfield $K$) which is determined by formulas $\{f(x)\neq 0: f(X)\in K\{X\}\}$.



The theory of differentially closed fields is model theoretically very nice. Namely, it admits elimination of quantifiers, elimination of imaginaries, it is complete and model complete. Further, $\DCF_0$ is $\omega$-stable with Morley rank $\omega$. Every differential field $K$ has a \textit{differential closure} which is defined as the prime model of $\DCF_0$ over $K$. The prime model always exists and is unique up to isomorphism (over $K$) in $\omega$-stable theories. We will denote \textit{the} differential closure of $K$ by $K^{\dif}$, while $K^{\alg}$ will denote the field theoretic algebraic closure. 



Furthermore, it is easy to see that the Morley degree of $\DCF_0$ is $1$. This means
that in a model $\mathcal{F} \models \DCF_0$ any definable set $A\subseteq F$ is either of finite rank
(it is \textit{small}) or it has rank $\omega$ and its complement has finite rank (it is \textit{big}).

Now we define differential curves and make some easy observations about them that will be used later in the paper.

\begin{definition}
A \textit{differential algebraic curve} $E$ in a differential field $K$ is a set in $K^2$ defined by a differential equation of two variables, i.e. $E = \{ (x,y)\in K^2 : f(x,y)=0\}$ for some $f(X,Y)\in K\{X,Y\}$. For brevity we will sometimes say \textit{differential curve} instead of differential algebraic curve.
\end{definition}

Note also that by an algebraic curve we mean a set defined by an algebraic equation of two variables. Let $\mathbb{D} := \{ (x,\D x) : x \in K \}$ be the graph of $\D$ in $K$. This is an example of a differential curve.



\begin{definition}\label{curve-general}
A \textit{differential curve in general sense} in a differentially closed field $\cal{F}$ is a definable subset of $F^2$ the generic fibres of which are of finite Morley rank. 
\end{definition}

Clearly any proper differential curve is a curve in general sense. On the other hand it is easy to notice that any curve $E$ in general sense must be contained in a proper differential curve. This means it must be defined by a formula of the form $\varphi(x,y) = [f(x,y)=0 \wedge \psi(x,y)]$ where $f$ is a differential polynomial and $\psi$ is any formula. Indeed, otherwise $E$ will contain a set of the form $f(x,y)\neq 0$ the generic fibres of which have rank $\omega$.  

{We could alternatively define curves in general sense to be definable sets (in $F^2$) of Morley rank less than $\omega \cdot 2$. The above argument shows that this is equivalent to the above definition. Thus, if $(a,b)$ is a pair of differentially independent elements and $\neg\varphi(a,b)$ holds in $\cal{F}$ then $\varphi(x,y)$ defines a curve in general sense.\footnote{One can also require $\MR(E)$ to be at least $\omega$ in order to avoid any degeneracies like $\D x = 0 \wedge \D y=0$ (which correspond to finite sets in $\ACF_0$), but it is not important for us.}} 




Finally let us fix some notations. We will use upper-case letters $X,Y,\ldots$ with possible subscripts for indeterminates of polynomials.  We will use lower-case letters for elements of a set and for variables in formulas (it will be clear from the context which one we mean). In particular if $f(X) \in F\{ X \}$ is a differential polynomial then $f(X)=0$ means that $f$ is identically zero, while $f(x)=0$  means $f$ vanishes at $x$ (or it is a formula with a free variable $x$).

\section{Definable derivations}\label{DefinableDerivations}
\setcounter{equation}{0}

{If $\D$ is a derivation on a field $(F;+,\cdot,0,1)$ then for any element $a\in F$ the map $a\cdot \D$ will be a derivation as well. We show in this section that in a differentially closed field all definable derivations are of that form.}
{
\begin{theorem}\label{def-deriv}
Let $\mathcal{F}=(F;+,\cdot,0,1,\D)$ be a differentially closed field and $\tilde{\D}$ be a definable (possibly with parameters) derivation. Then there exists an element $a\in F$ such that $\tilde{\D}=a\D$.
\end{theorem}}

{Though this fact is well known (a proof can be found for example in \cite{su}), we nevertheless present our proof here as we are going to need it in Section \ref{Further examples}.}

{The following well-known result is a characterisation of definable functions in a differentially closed field (see, for example, \cite{Pillay1} or \cite{TZ}, Exercise 6.1.14).}

{
\begin{lemma}\label{L3-14}
Let $\mathcal{F}$ be a differentially closed field and $f:F^k\rightarrow F$ be a definable (possibly with parameters) function in $\mathcal{F}$. Then there is a partition of $F^k$ into a finite number of definable subsets $U_i$ such that $f$ is given by a differential rational function on each of them (this means, in particular, that each of these rational functions is determined on the corresponding set).
\end{lemma}}

{We establish one more result before proving the Theorem \ref{def-deriv}.
\begin{lemma}\label{L3-15}
Suppose $\D$ and $\D_1$ are derivations on a field $(F;+,\cdot,0,1)$ such that there is $t\in F$ with $\D t=1$. Let $P(X_0,\ldots,X_n,Y)$  be a non-zero polynomial over $F$ such that
\begin{equation}
P(X,\D X,\ldots,\D^nX,\D_1X)= 0.
\end{equation}
Then $\D_1=a\cdot \D$, where $a=\D_1t$.
\end{lemma}}
\begin{proof}{
For an element $x\in F$ and an arbitrary rational number $r$ one has $P(x+r,\D x,\ldots,\D^nx,\D_1x)= 0,$ hence
$$P(X,\D x,\ldots,\D^nx,\D_1x)= 0$$
(as a polynomial of $X$). Therefore all coefficients of this polynomial are zeros. Since $P(X_0,\ldots,X_n,Y)$ is non-zero, if we consider it as a polynomial of $X_0$, it will have a non-zero coefficient that is a polynomial of $X_1,\ldots,X_n,Y$. It must vanish at $(\D x,\ldots,\D^nx,\D_1x)$. This is true for all $x\in F$.}

{Thus for a non-zero polynomial $P_1$ we have
\begin{equation*}
P_1(\D X,\ldots,\D^nX,\D_1X)= 0.
\end{equation*}
Again, fixing an element $x\in F$ we see that for any rational $r$ one has
$P_1(\D x+r,\D^2x,\ldots,\D^nx,\D_1x+ar)= 0$ (we substitute $X=x+rt$). This implies
$$P_1(X,\D^2x,\ldots,\D^nx,\D_1x-a\D x+aX)= 0.$$
Replacing $X$ by a fixed element $y\in F$ and taking $x+rt^2$ instead of $x$ we get
$$P_1(y,\D^2x+2r,\D^3x,\ldots,\D^nx,\D_1x-a\D x+ay)= 0.$$
Therefore $$P_1(y,X,\D^3x,\ldots,\D^nx,\D_1x-a\D x+ay)= 0.$$
Arguing as above we show that for some non-zero polynomial $P_2$ we have
$$P_2(y,\D^3x,\ldots,\D^nx,\D_1x-a\D x+ay)= 0$$
for all $x,y\in F$.
Proceeding this way one can prove that there is a non-zero polynomial $Q(Z_1,Z_2)\in F[Z_1,Z_2]$ such that
$$Q(Y,\D_1X-a\D X+aY)= 0.$$
Now suppose for some $u\in F$ we have $\D_1u\neq a\D u$. Then for any natural number $n$ one has $\D_1(nu)\neq a\D(nu)$. This means that for any $y\in F$ the polynomial $Q(y,ay+Z)$ equals zero for infinitely many values of $Z$, hence, it is identically zero. This yields $Q(Y,Z)= 0$. We arrived at a contradiction, therefore $\D=a\D_1$.}
\end{proof}

\begin{proof}[Proof of Theorem \ref{def-deriv}]


{
From Lemma \ref{L3-14} it follows that there are definable sets $U_i\subseteq F$ such that $\tilde{\D}$ is given by a differential rational function on each $U_i$. Therefore there are differential polynomials $f_i(X),g_i(X)\in F\{X\}$ such that $f_i(x)\cdot \tilde{D}(x) = g_i(x)$ and $f_i(x)\neq 0$ for all $x\in U_i$. We know that $f_i(X)=P_i(X,\D X,\ldots,\D^mX),\ g_i(X)=Q_i(X,\D X,\ldots,\D^mX)$ for some polynomials $P_i$ and $Q_i$ over $F$. Form the polynomial
$$P(X_0,\ldots,X_m,Y)=\prod_i (P_i(X_0,\ldots,X_m)\cdot Y - Q_i(X_0,\ldots,X_m)).$$
This is a non-zero polynomial and
$$P(X,\D X,\ldots,\D^mX,\tilde{\D}X) = 0.$$
As $\mathcal{F}$ is differentially closed, there exists $t\in F$ with $\D t=1$. Now Lemma \ref{L3-15} yields the desired result.}
\end{proof}


\section{Model theoretic properties of the reducts}\label{ModelTheoreticProperties}
\setcounter{equation}{0}

From now on we will work in a differentially closed field $\mathcal{F} = ( F;+,\cdot,0,1,\D )$ which we will assume to be sufficiently saturated. Thus, it will serve as a monster model for us.

For a collection $R$ of definable sets in (Cartesian powers of) $\mathcal{F}$,  we define the $R$-\textit{reduct} $\mathcal{F}_{R}$ of $\mathcal{F}$ to be the structure $(F;+,\cdot,0,1,P)_{P \in R}$ in the language $\mathfrak{L}_{R} = \{ +,\cdot,0,1\} \cup R$ (the elements of $R$ are relation symbols in the language $\mathfrak{L}_{R}$). We will omit $R$ and just say ``reduct'' whenever no confusion can arise. We will say that $R$ (or the reduct $\mathcal{F}_{R}$) is \textit{algebraic} if all relations of $R$ can be defined in the pure field $(F;+,\cdot,0,1)$. If $R$ consists of just one relation $E$ then we will write $\mathcal{F}_{E}$ for the corresponding $E$-reduct.

In this section we examine basic model theoretic properties of the reducts $\mathcal{F}_{R}$. Though we will sometimes assume $R$ is finite, most of our results will be valid for an arbitrary $R$. From the point of view of Question \ref{Question-definability} the assumption of finiteness of $R$ is no loss of generality as a possible definition of $\D$ would anyway contain only finitely many occurrences of relation symbols from $R$. 


We start by introducing a piece of notation. In order to distinguish between the same concepts in the differentially closed field $\mathcal{F}$ and in the reduct $\mathcal{F}_{R}$, we will add a subscript $\D$ or $R$ respectively to their notations. Thus $\MR_{\D},~ \MD_{\D},~ \tp_{\D},$ $\dcl_{\D},~ \acl_{\D}$ stand for Morley rank, Morley degree, type, definable closure and algebraic closure respectively in $\mathcal{F}$ while $\MR_{R},~ \MD_{R},~ \tp_{R},~ \dcl_{R},~ \acl_{R}$ stand for the same notions in $\mathcal{F}_{R}$. 

Also we will need to consider generic elements and types. \textbf{By generic we will always mean generic in the differentially closed field $\mathcal{F}$ (rather than in $\mathcal{F}_R$) unless explicitly stated otherwise.} If we do not specify over which set an element is generic then we mean over the empty set.\\

Finally, we turn to model theoretic properties of the reducts. Clearly $\mathcal{F}_{R}$ is an $\omega$-stable structure. We now find its Morley rank.

\begin{proposition}\label{rank}
$\mathcal{F}_{R}$ has Morley rank $\omega$ unless $R$ is algebraic.
\end{proposition}
\begin{proof}
First of all, since $\mathcal{F}_R$ is a reduct of $\mathcal{F}$, and the latter has Morley rank $\omega$, we have $\MR(\mathcal{F}_R) \leq \omega$. So we need to prove $\MR(\mathcal{F}_R) \geq \omega$.

It suffices to prove this for $R=\{P\}$ where $P$ is a non-algebraic unary relation which has finite Morley rank in the differentially closed field $\mathcal{F}$. The case $P=C$ (the field of constants) is a well known example. In this case the reduct is just an algebraically closed field with a unary predicate for an algebraically closed subfield. Our proof below is an adaptation of a known proof for this special case (see, for example, \cite{Mar}, exercise 6.6.17, d).

As $P$ is non-algebraic, it must be infinite and hence $\MR_{R}(P) \geq 1$. Also $P$ has finite Morley rank in $\mathcal{F}$, so $(\mathbb{Q}(P))^{\alg}\neq F$. Now for an element $x \in F \setminus (\mathbb{Q}(P))^{\alg}$ define
$$X_n = \left\{ y \in F : \exists a_0,\ldots,a_{n-1} \in P~ \left(y = \sum a_i x^i \right) \right\}.$$
The map $\pi : P^{n+1} \rightarrow X_{n+1}$ given by $(a_0,\ldots,a_n) \mapsto a_0 + a_1x+ \ldots +a_n x^n $ is a definable bijection. Hence $\MR_{R}(X_n) = \MR_{R} (P^n) \geq n$. Therefore $\MR_{R}(F)=\omega$.
\end{proof}

We will assume throughout the paper that $R$ is not algebraic and so $\mathcal{F}_R$ has Morley rank $\omega$.

\begin{remark}\label{High-rank-formulas}
As we saw in the proof, if $a \in F$ is a differentially transcendental element then for each $n < \omega$ there is a definable (in $\mathcal{F}_R$) set $X_n \subseteq F$, defined over $a$, such that $n \leq \MR_R(X_n) < \omega$.
\end{remark} 

Further, observe that $\mathcal{F}_{R}$ has Morley degree $1$. If $\varphi(x)$ is a formula {(of one variable)} in the language $\mathfrak{L}_{R}=\{+,\cdot,0,1\} \cup R$ then in the language $\mathfrak{L}_{\D}$ it is equivalent to a quantifier-free formula. If it is an equation in conjunction with something else, then $\MR_{R} (\varphi) < \omega$ otherwise $\MR_{R} (\varphi) = \omega$. Also, $\MR_{R} (\varphi) \leq \MR_{\D} (\varphi)$ and these ranks are finite or infinite simultaneously. {Indeed, if $\MR_{\D} (\varphi) = \omega$ then $\MR_{\D} (\neg \varphi ) < \omega$, and so $\MR_R(\neg \varphi) < \omega$. Therefore, $\MR_R(\varphi) = \omega$ since $\MR_R(x=x)=\omega$ as proven above.}

There is a unique generic $1$-type in $\mathcal{F}_{R}$ given by
$$ \{ \varphi(x) : \varphi \in \mathfrak{L}_{R},~ \MR_{R} (\varphi) = \omega \} = \{ \neg\varphi(x) : \varphi \in \mathfrak{L}_{R},~  \MR_{R} (\varphi) < \omega \}.$$

Similarly, the unique generic $n$-type is given by formulas of Morley rank $\omega \cdot n$.

Now let us discuss the issue of quantifier elimination for $\mathcal{F}_{R}$. First notice that even when $R = \{ \mathbb{D} \}$, $\mathcal{F}_{R}$ does not admit quantifier elimination for $y=\D^2x$ is existentially definable but not quantifier-free definable. It turns out that this is a general phenomenon.

\begin{corollary}\label{noqf}
If $R$ is non-algebraic and finite then the reduct $\mathcal{F}_{R}$ does not admit elimination of quantifiers.
\end{corollary}
\begin{proof}
Suppose $R$ is not algebraic but $\mathcal{F}_{R}$ has quantifier elimination. Then any formula with one free variable must be equivalent to a Boolean combination of algebraic polynomial equations (in the language of rings) and formulas of the form
$$Q(p_1(x),\ldots,p_n(x))$$
where $Q \in R$ is an $n$-ary predicate and $p_i$'s are algebraic polynomials.
But clearly if such a formula has finite Morley rank then the latter is uniformly bounded, i.e. there is a bound which is the same for all formulas of finite Morley rank (remember that $R$ is finite). This contradicts Proposition \ref{rank}.
\end{proof}

One sees that although in the case $R=\{ \mathbb{D} \}$ the reduct does not have quantifier elimination, it is nevertheless model complete. In general it is true if $\D$ is existentially definable. We show this below.

\begin{lemma}\label{exist-univ}
Let $\mathcal{M}$ be a structure. If a function $f:M^n \rightarrow M$ is existentially definable in $\mathcal{M}$ then it is also universally definable.
\end{lemma}
\begin{proof}
If $\phi(\bar{x},y)$ defines $f$ then so does $\forall z (z=y \vee \neg \phi(\bar{x},z))$.
\end{proof}

\begin{proposition}\label{modelcomplete}
If $\D$ is existentially definable in $\mathcal{F}_{R}$ then $T_{R}:= \Th (\mathcal{F}_R)$ is model complete.
\end{proposition}
\begin{proof}
Suppose that $\D$ is existentially definable. Take an arbitrary formula $\varphi \in \mathfrak{L}_{R}$. In the language of differential rings it is equivalent to a quantifier-free formula, i.e. to a Boolean combination of differential equations. Each differential equation is existentially definable in the reduct and, by Lemma \ref{exist-univ}, it is also universally definable. Substituting existential definitions in positive parts (i.e. equations) and universal definitions in negative parts (inequations), we get an existential formula in the language $\mathfrak{L}_{R}$. Thus any formula in the language of the reduct is equivalent to an existential formula. This is equivalent to model completeness.
\end{proof}

Thus, model completeness is the deepest possible level of quantifier elimination that we can have for $T_R$. As we will see in the last section, under a natural assumption, definability of $\D$ will imply that $T_R$ is model complete.

\section{An example}\label{some examples}
\setcounter{equation}{0}

In this section we show that in a certain class of reducts $\D$ is definable. It will be used later to establish some criteria for definability of $\D$.

Choose an element $t\in F$ with $\D t=1$ (it exists because our field is differentially closed) and add it as a constant symbol to our language. Thus from now on we work in the language $\{+,\cdot,\D,0,1,t\}$ for differential fields, which by abuse of notation we will again denote by $\mathfrak{L}_{\D}$. Correspondingly all reducts will be considered in the language $\mathfrak{L}_{R}=\{+,\cdot,0,1,t\} \cup R$. Again abusing the nomenclatures we will call $\mathfrak{L}_{\D}$ the language of differential rings and $\mathfrak{L}_{R}$ the language of the reducts. This means that we do not count $t$ as a parameter in our formulas, i.e. we are free to use $t$ in formulas and declare that something is definable without parameters. Note that this does not affect any of the results proved in the previous section. Let us also mention that after adding $t$ to our language (and requiring that a derivation takes the value $1$ at $t$) the only candidate for a definable derivation can be $\D$ (see Theorem \ref{def-deriv}).

For a formula $\varphi(\bar{x})$ in the language $\mathfrak{L}_{R,\D} = \mathfrak{L}_{R} \cup \mathfrak{L}_{\D}$ and a tuple $\bar{a} \in F$ we will sometimes write $\mathcal{F} \models \varphi(\bar{a})$. This is an abuse since in general $\varphi$ is not in the language of differential rings, but clearly $\mathcal{F}$ can be canonically made into an $\mathfrak{L}_{R,\D}$-structure. 

In general, if the relations in $R$ are defined with parameters and $\D$ is definable then it will be definable with parameters as well. But in many cases we do not use any extra parameters to define $\D$. So for simplicity we will assume that $R$ consists of $0$-definable relations in $\mathcal{F}$, i.e. relations defined over $k_0=\mathbb{Q}(t)=\dcl(\emptyset)$. Thus from now on by definable we will mean definable without parameters unless explicitly stated otherwise.

We denote the theory of the reduct by $T_{R} := \Th (\mathcal{F}_{R})$. We will say that there is a derivation $\D_K$ on a model $\mathcal{K}_{R} \models T_{R}$ which is \textit{compatible} with $R$. This means that $(K;+,\cdot,\D_K, 0,1,t, P)_{P \in R} \equiv (F;+,\cdot,\D, 0,1,t, P)_{P \in R}$, i.e. the differential field $\mathcal{K}=(K;+,\cdot,\D_K, 0,1,t)$ is differentially closed with $\D_Kt=1$ and the sets from $R$ are defined by the same formulas as in $\mathcal{F}$.

Throughout the paper we let $E$ be a differential curve (possibly in general sense); as we noted above the corresponding reduct will be denoted $\mathcal{F}_{E}$. Recall also that $\mathbb{D}=\{ (x, \D x): x \in F \}$ is the graph of $\D$.

Now we prove an auxiliary result which will be used several times throughout the paper. {It states that $(\mathbb{Q}(t))^n$ is Kolchin-dense in $F^n$ for each $n$.}

\begin{lemma}[{cf. \cite{Mar-dif}, Lemma A.4}]\label{nonzero}
For any non-zero differential polynomial $f(X_1,\ldots,X_n)$ over $\mathbb{Q}(t)$ there are elements $t_1,\ldots,t_n \in \mathbb{Q}[t]$ such that $f(t_1,\ldots,t_n) \neq 0$.
\end{lemma}
\begin{proof}
{First assume $f$ is a polynomial of one variable $X$.} Let $\ord(f)=n$. Since $\mathcal{F}$ is differentially closed, we can find an element $u \in F$ with $\D^{n+1} u =0 \wedge f(u) \neq 0$. Then clearly 
$$u=c_n t^n + \ldots + c_1t +c_0$$
for some constants $c_0,\ldots,c_n \in C$.

Now for constants $\lambda_0, \ldots, \lambda_n$ denote 
$$p(t, \bar{\lambda})=\lambda_n t^n + \ldots + \lambda_1t +\lambda_0.$$
Since $t$ is transcendental over $C$,  there are algebraic polynomials $q_i(X_0, \ldots, X_n) \in \mathbb{Q}[X_0,\ldots,X_n],~ i=1,\ldots,m,$ such that for all $\bar{\lambda} \in C^{n+1}$
$$f(p(t, \bar{\lambda}))=0 \mbox{ iff } \bigwedge_{i=1}^m q_i(\bar{\lambda})=0.$$
Let $V \subseteq C^{n+1}$ be the algebraic variety over $\mathbb{Q}$ defined by $\bigwedge_{i=1}^m q_i(\bar{\lambda})=0$. Then as we saw above $V(C) \neq C^{n+1}$, and hence $V(\mathbb{Q}) \subsetneq \mathbb{Q}^{n+1}$. So there is a tuple $\bar{r} \in \mathbb{Q}^{n+1}$ with $\bar{r} \notin V(\mathbb{Q})$. Therefore $f(p(t,\bar{r}))\neq 0$  and $p(t,\bar{r}) \in \mathbb{Q}[t]$.

{Now we prove the general case (when $f$ has more than one variables) by induction on $n$. If $f = f(X_1,\ldots,X_n)$ with $n>1$ then consider it as a differential polynomial $g(X_1,\ldots,X_{n-1})$ of $n-1$ variables over the differential ring $\mathbb{Q}(t)\{ X_n \}$. Choose a non-zero coefficient of $g$ which will be a non-zero differential polynomial $h(X_n) \in \mathbb{Q}(t) \{ X_n \}$. As we proved above there is $t_n \in \mathbb{Q}[t]$ such that $h(t_n) \neq 0$. Now the polynomial $f(X_1,\ldots,X_{n-1},t_n)$ is a non-zero polynomial of $n-1$ variables over $\mathbb{Q}(t)$ and we are done by the induction hypothesis.}
\end{proof}

\begin{remark}
The proof shows that we can choose $t_1,\ldots,t_n$ from $\mathbb{Z}[t]$ (and even from $\mathbb{N}[t]$).
\end{remark}

\begin{definition}\label{multi-degree}
Introduce the reverse lexicographical order on $(n+1)$-tuples of integers, that is, $(\alpha_0,\ldots,\alpha_n)< (\beta_0,\ldots,\beta_n)$ if and only if for some $j$, $\alpha_i = \beta_i$ for $i> j$ and $\alpha_j<\beta_j$. The \textit{multi-degree} of an algebraic polynomial $Q(X_0,\ldots,X_n)$ is   the greatest (with respect to this order) $(n+1)$-tuple $(\alpha_0,\ldots,\alpha_n)$  for which $X_0^{\alpha_0}\cdot \ldots \cdot X_n^{\alpha_n}$ appears in $Q$ with a non-zero coefficient. The multi-degree of a differential polynomial $f(X)=P(X,\D X,\ldots,\D^nX)$ is defined as that of $P$. 
\end{definition}

\begin{theorem}\label{curve-def}
If $E$ (a differential algebraic curve) contains the graph of $\D$ then $\D$ is quantifier-free definable in $\mathcal{F}_E$.
\end{theorem}

Before proving the theorem we give an example which helps to understand how the proof works.

\begin{example}
Suppose $E$ is given by $(\D y-\D^2x) \cdot \D x=0$. Then $E(x+t,y+1)=[(\D y-\D^2x) \cdot (\D x+1)=0]$. The conjunction $E(x,y)\wedge E(x+t,y+1)$ implies $\D y-\D^2x=0$. Now we substitute $ x\mapsto tx,~ y\mapsto x+ty$ and get $t(\D y-\D^2 x)+y-\D x=0$. Subtracting the previous equation multiplied by $t$ we get $y-\D x =0$. Thus the formula $E(x,y)\wedge E(x+t,y+1)\wedge E(tx,x+ty)\wedge E(tx+t,x+ty+1)$ defines $\D$.
\end{example}

\begin{proof}\footnote{I am grateful to Ehud Hrushovski for detecting a gap in the initial version of the proof and helping me to fix it.}
{Let $E$ be given by a differential equation $f(x,y)=0$. We know that $f(X,\D X)$ identically vanishes. Denote $U:= Y-\D X$ and consider the differential polynomial $g(X,U):=f(X,U+\D X)$. Clearly $g(X,0)=0$. }

{First we intersect additive translates to ``eliminate'' $x$ and define a differential equation $h(u)=0$ for some differential polynomial $h(U)$. If $g(X,U)$ depends on $X$ (i.e. $g(X,U) \in k_0\{ X,U\} \setminus k_0\{ U\}$) then we can find (see Lemma \ref{nonzero}) $p(t)\in \mathbb{Q}[t]$ such that $g(X+p(t),U)\neq g(X,U)$. Clearly, $U$ is invariant under the transformation $X\mapsto X+p(t),~ Y\mapsto Y+p'(t)$ where $p'(Z)=\frac{\partial p}{\partial Z}$. So consider the formula $E(x,y)\wedge E(x+p(t),y+p'(t))$. It is equivalent to $g(x,u)=0 \wedge g(x+p(t),u)=0$ which implies $g_1(x,u):=g(x,u)-g(x+p(t),u)=0$. The leading terms of the differential polynomials $g(X,U)$ and $g(X+p(t),U)$ in variable $X$ (i.e. the sums of monomials in these polynomials that have highest multi-degree in $X$) are the same and hence they cancel out in the difference $g_1(X,U):=g(X,U) - g(X+p(t),U)$. On the other hand $g_1(X,U)\neq 0$ by our choice of $p$ and the multi-degree of $g_1$ in $X$ is strictly less than that of $g$. In other words, if the multi-degree of $g$ in $X$ is bigger than $(0,\ldots,0)$ then we can reduce it. Now if $g_1(X,U)$ depends on $X$ then we do the same for $g_1$. We keep repeating this process and reduce the multi-degree of our differential polynomial step by step until it becomes $(0,\ldots,0)$. This means we get a curve $h(u)=0$ for a non-zero differential polynomial $h$, which contains a quantifier-free definable set in our reduct. It is also clear that the latter contains the curve $u=0$ (the graph of $\D$).}

{Now we use multiplicative translates to define the curve $u=0$ (which is actually $y=\D x$).
Let $p(t)\in \mathbb{Q}[t]$. When we substitute $X \mapsto p(t)X,~ Y\mapsto p'(t)X+p(t)Y$ then $U$ is replaced by $p(t)U$. Then $h(u)=0 \wedge h(p(t)u)=0$ is implied by a quantifier-free formula in the language of the reduct and implies $h_{\alpha,1}(u):=p(t)^{\alpha} h(u) - h(p(t)u) = 0$ for any positive integer $\alpha$. If $(\alpha_0,\ldots,\alpha_n)$ is the multi-degree of $h$ then taking $\alpha:=\alpha_0+\ldots+\alpha_n$ the leading terms of the differential polynomials $p(t)^{\alpha} h(U)$ and $h(p(t)U)$ will coincide and will cancel out in the difference $h_{\alpha,1}(U):=p(t)^{\alpha} h(U) - h(p(t)U)$. By an appropriate choice of $p$ we can also guarantee that $h_{\alpha,1}(U)$ is non-zero unless $h(U)=h(1)\cdot U^{\alpha}$. Indeed, if $h(U)\neq h(1)\cdot U^{\alpha}$ then the polynomial $h(V\cdot U)-V^{\alpha}\cdot h(U)$ is non-zero and hence there is $p(t) \in \mathbb{Q}[t]$ such that $h(p(t)\cdot U) \neq p(t)^{\alpha}\cdot h(U)$, therefore $h_{\alpha,1}(U)\neq 0$. Thus, if $h(U)$ is not a homogeneous algebraic polynomial then $h_{\alpha,1}$ is non-zero and its multi-degree is strictly less than that of $h$. Now if $h_{\alpha,1}(U)$ is not algebraic homogeneous then we repeat the above procedure for $h_{\alpha,1}$. Iterating this process we will eventually obtain an equation $u^{\alpha}=0$ for some positive integer $\alpha$ which is equivalent to $u=0$. Taking into account that all the sets defined this way contain $u=0$ we see that at the last step we have defined $u=0$ which, in terms of $x$ and $y$, is the curve $y=\D x$.}

Finally note that we only take conjunctions of atomic formulas here, hence the definition is quantifier-free.
\end{proof}

\begin{remark}
Strictly speaking, for the ``quantifier-free'' part of the theorem to be true we need to pick $p(t)\in \mathbb{N}[t]$ each time. Alternatively, we could add unary functions for multiplicative and additive inverses to our language.
\end{remark}


\begin{corollary}\label{general-curve}
If $E$ is a curve in general sense that contains $\mathbb{D}$ then $\D$ is quantifier-free definable.
\end{corollary}
\begin{proof}
Being a curve in general sense, $E$ is defined by a formula of the form $f(x,y)=0 \wedge \psi(x,y)$ for $\psi$ a quantifier free formula in the language of differential fields. 
Now for the curve $E'$ given by the equation $f(x,y)=0$ we have a definition of $\D$. Suppose it is given by the formula $\varphi(x,y)$ in the reduct $\mathcal{F}_{E'}$. We claim that the same formula defines $\D$ in $\mathcal{F}_{E}$. Indeed, as we take only conjunctions to define $\D$ from $E'$, the set defined by $\varphi(x,y)$ in $\mathcal{F}_{E}$ will be contained in $\mathbb{D}$. On the other hand it clearly contains $\mathbb{D}$. Therefore it defines $\mathbb{D}$.
\end{proof}

We will give further examples and non-examples (of differential equations defining $\D$) in Section \ref{Further examples}, but first we need to establish some facts on generic points which we do in the next section.

\section{Generic points}\label{generic points}
\setcounter{equation}{0}

Recall that we work in a saturated differentially closed field $\mathcal{F}$. From now on we fix a generic (in the sense of $\DCF_0$, that is, differentially transcendental) point $a\in F$. We first prove that if $\D a$ can be defined from $a$ then we can recover the whole of $\D$.


\begin{proposition}\label{generic}
Suppose a formula $\varphi(x,y) \in \mathfrak{L}_{R}$ defines $\D a$ from $a$, that is, 
$$\mathcal{F} \models \forall y ( \varphi(a,y) \leftrightarrow y=\D a).$$ 
Then $\D$ is definable (without parameters). Moreover, if $\varphi$ is existential then $\D$ is existentially definable.
\end{proposition}
\begin{proof}[First proof]
First of all observe that since the generic type is unique, for any differentially transcendental element $b\in F$ we have
$$\mathcal{F} \models \forall y ( \varphi(b,y) \leftrightarrow y=\D b).$$ 

Let $A$ be the set defined by $\varphi(x,y)$ and define
$$B:=\{ (b,\D b) : b \mbox{ generic in } \mathcal{F} \} \subseteq A.$$

At generic points $b$ the formula $\varphi$ defines $\D b$ but we do not have any information about non-generic points. So we need to shrink the set $A$ to a subset of $\mathbb{D}$ in order to avoid any possible problems at non-generic points. The set $A$, being a curve in general sense {(its fibre over any generic point $x=b$ consists of one element and hence is small)}, must be defined by a formula $f(x,y)=0 \wedge \psi(x,y)$ (in the language of differential rings). Then $f(a,\D a )=0$ and hence $f(X,\D X) = 0$. Therefore $\D$ can be defined from the differential curve $f(x,y)=0$ by Theorem \ref{curve-def}. Taking into account that for a generic element $b$ the elements $b+p(t)$ and $p(t)b$ are generic as well for any $p(t)\in \mathbb{Q}[t]\setminus \{ 0\}$, we see that the sets $\varphi(x,y) \wedge \varphi(x+p(t),y+p'(t))$ and $\varphi(x,y) \wedge \varphi(p(t)x,p(t)y+p'(t)x)$ contain $B$. Arguing as in the proofs of Theorem \ref{curve-def} and Corollary \ref{general-curve}, after taking sufficiently many conjunctions of such formulas we will eventually define a set $B'$ such that it contains $B$ and is contained in the graph $\mathbb{D}$ of $\D$. Note that $B'$ is $0$-definable.

Treating $\mathbb{D}$ as an additive group we prove the following.
\begin{claim}
$\mathbb{D}=B'+B'$.
\end{claim}

Clearly $B'+B' \subseteq \mathbb{D}$. Let us show that the converse inclusion holds. Any element $d \in F$ has a representation $d=b_1+b_2$ with $b_1$ and $b_2$ generic. Indeed, take $b_1$ to be generic over $d$ and choose $b_2=d-b_1$. Hence $(d, \D d) = (b_1, \D b_1) + (b_2, \D b_2) \in B+B \subseteq B'+B'$.

This gives a definition of $\D$ without parameters. Moreover, if $\varphi$ is existential then we get an existential definition.
\end{proof}

\begin{remark}
The group $\mathbb{D}$ is in fact a connected $\omega$-stable group (its Morley degree is one). Therefore the equality $\mathbb{D}=G+G$ holds for any definable subset $G$ of $\mathbb{D}$ with $\MR(G)=\MR(\mathbb{D})$ (see, for example, \cite{Mar}, Chapter 7, Corollary 7.2.7). We could use this to show that $\mathbb{D}=B'+B'$ since $\MR(B')=\MR(\mathbb{D})=\omega$. In fact, the idea is the same as in the above claim; one just passes to a saturated extension and uses the above argument there. 
 
\end{remark}

We will shortly give another proof to Proposition \ref{generic}. For this we first observe that if $\D$ is definable with independent parameters then it is also definable without parameters.

\begin{lemma}\label{generic-withoutparameters}
Suppose $\psi(x,y,u_1,\ldots,u_n) \in \mathfrak{L}_{R}$ and $b_1,\ldots,b_n$ are differentially independent elements in $\mathcal{F}$. If the formula $\psi(x,y,\bar{b})$ defines $y=\D x$ then there are $0$-definable elements $t_1,\ldots,t_n \in k_0=\mathbb{Q}(t)$ such that $\psi(x,y,\bar{t})$ defines $\D$ (and so $\D$ is $0$-definable).
\end{lemma}
\begin{proof}
We have
$$\mathcal{F} \models \psi(x,y,\bar{b}) \longleftrightarrow y=\D x.$$
Therefore
$$q(\bar{z}):=\tp_{\D}(\bar{b}) \models \psi(x,y,\bar{z}) \longleftrightarrow y=\D x.$$
Since $q(\bar{z})$ is the generic $m$-type in $\DCF_0$, it consists only of differential inequations. Applying compactness and taking into account that conjunction of finitely many inequations is an inequation as well, we conclude that there is a differential polynomial $f(Z_1,\ldots,Z_m)$ over $k_0$ such that
$$\mathcal{F} \models \forall \bar{z} (f(\bar{z}) \neq 0 \longrightarrow \forall x,y(\psi(x,y,\bar{z}) \leftrightarrow y=\D x)).$$

By Lemma \ref{nonzero} we can find elements $t_1,\ldots,t_m \in k_0$ such that $f(t_1,\ldots,t_m)$ is non-zero. Now we see that
$$\mathcal{F} \models \psi(x,y,\bar{t}) \longleftrightarrow y=\D x$$
and we are done.
\end{proof}

\begin{proof}[Second proof of Proposition \ref{generic}]
{Let $(b_1,b_2) \in F^2$ be a differentially independent tuple. Then for every $d\in F$ the differential transcendence degree of $d,d+b_1,d+b_2$ is at least $2$.} It is easy to deduce from this that the following formula defines $\D$:
\begin{gather*}
\exists u_1,u_2 (\varphi(b_1,u_1) \wedge \varphi(b_2,u_2) \wedge [ (\varphi(x,y)\wedge \varphi(x+b_1,y+u_1)) \\ \vee (\varphi(x,y)\wedge \varphi(x+b_2,y+u_2))\vee (\varphi(x+b_2,y+u_2)\wedge \varphi(x+b_1,y+u_1)) ] ).
\end{gather*}

Now Lemma \ref{generic-withoutparameters} concludes the proof.
\end{proof}

The idea that the behaviour of $\D$ at generic (differentially transcendental) points determines its global behaviour as a function can be developed further. We proceed towards this goal in the rest of this section. 

Next we show that if $\D a$ is not generic over $a$ (in the reduct) then it is in fact definable and hence $\D$ is definable. Let $p(y):=\tp_{R}(\D a/a)$ be the type of $\D a$ over $a$ in $\mathcal{F}_R$.

\begin{theorem}\label{finiterank}
The derivation $\D$ is definable in $\mathcal{F}_{R}$ if and only if $p$ has finite Morley rank (in $\mathcal{F}_R$).
\end{theorem}
\begin{proof}
Obviously, if $\D$ is definable then $p$ is algebraic and hence has Morley rank $0$. Let us prove the other direction.

Let $\varphi(a,y)\in p$ be a formula of finite Morley rank. Trivially $\mathcal{F} \models \varphi(a,\D a)$ and  $\varphi(x,y)$ defines a curve in general sense.
As in the proof of Proposition \ref{generic} we can define a big subset $\psi(x,y)$ of $\mathbb{D}$, that is, a subset of Morley rank $\omega$. This set certainly contains the point $(a,\D a)$ and $\psi(a,y)$ defines $\D a$. Thus $\D a$ is definable over $a$ and Proposition \ref{generic} finishes the proof.
\end{proof}

\begin{remark}\label{remark-existential-def}
The proof shows that if $\varphi(x,y)$ is an existential formula of rank $< \omega \cdot 2$ which is true of $(a,\D a)$ then $\D$ is existentially definable.
\end{remark}

\begin{corollary}
In the reduct, $\D a$ is either generic or algebraic (in fact, definable) over $a$.
\end{corollary}

\begin{lemma}
If $p$ is isolated then it has finite Morley rank (in the reduct).
\end{lemma}
\begin{proof}
The argument here is an adaptation of the proof of the fact that in differentially closed fields the generic type is not isolated.

Suppose $p$ is isolated but has rank $\omega$, i.e. it is the generic type over $a$ (in the reduct). Then
\begin{equation*}
p(y)=\{ \neg \varphi(a,y) : \varphi \in \mathfrak{L}_R,~ \mathcal{F} \models \varphi(a,\D a) \mbox{ and } \MR_R(\varphi(a,y))< \omega\}.
\end{equation*}
Suppose $\neg \psi(a,y)$ isolates $p$. By Remark \ref{High-rank-formulas} there is a formula $\varphi(a,y)$ for which $\MR_R(\psi(a,y)) < \MR_R(\varphi(a,y)) < \omega$. Then $\varphi(a,y) \wedge \neg \psi(a,y)$ is consistent. A realisation of this formula cannot be generic, for $\varphi$ has finite Morley rank. This is a contradiction.
\end{proof}

As an immediate consequence one gets the following result.

\begin{corollary}\label{isolated}
The derivation $\D$ is definable in $\mathcal{F}_{R}$ if and only if $p$ is isolated.
\end{corollary}

\begin{remark}
We can consider the quantifier-free type $q(y):=\qftp(\D a/a)$. Then $\D$ is quantifier-free definable if and only if this type is isolated, if and only if it has finite Morley rank. 
\end{remark}

Notice that in stability-theoretic language we have proved that $\D$ is definable if and only if $\tp_{R}(\D a/a)$ forks over the empty set. Indeed, $\MR_{R}(\D a) = \omega$ (since it is generic in the differentially closed field) and forking in $\omega$-stable theories means that Morley rank decreases, hence $\tp_{R}(\D a/a)$ forks over $\emptyset$ if and only if $\MR_{R}(\D a/ a) < \omega$. 
In terms of forking independence we have the following formulation: $\D$ is definable if and only if $a \nforkind_{} \D a$ in $\mathcal{F}_{R}$. This will be generalised in the next section. 
Note also that all the above results will remain true if we replace Morley rank everywhere with $\U$-rank.

Now add a differentially transcendental element $a$ to our language and consider the reducts in this new language. Denote the theory of $\mathcal{F}_{R}$ in this language by $T^+_{R}$. Assume that each model of $T^+_{R}$ comes from a differentially closed field, that is, each model $\mathcal{K}_{R}$ is the reduct of a differentially closed field $\mathcal{K}=(K;+,\cdot,\D_K,0,1,t,a)$ in which $a$ is generic (differentially transcendental) and relations from $R$ are interpreted canonically (i.e. they are defined in $\mathcal{K}$ by the same formulas as in $\mathcal{F}$). Then the type $p(y)$ will be realised by $\D_Ka$ in $\mathcal{K}_R$. The omitting types theorem now yields that $p$ must be isolated. Thus, we have established the following result.

\begin{theorem}\label{T+}
If each model of $T^+_{R}$ is the $R$-reduct (with canonical interpretation) of a model of $\DCF_0$, then $\D$ is definable.
\end{theorem}

In other words, this means that if each model of $T^+_{R}$ is equipped with a derivation which is compatible with $R$ then $\D$ is definable. The converse of this holds as well trivially.

This is similar to Beth's definability theorem in spirit (see \cite{Poizat}). Beth's theorem in this setting means that if each model of $T^+_{R}$ has at most one derivation compatible with $R$ then $\D$ is definable. We showed that if each model has at least one derivation then $\D$ is definable. Also it is worth mentioning that unlike Beth's definability theorem, this statement is not true in general for arbitrary theories. 

\section {Further examples}\label{Further examples}
\setcounter{equation}{0}
In this section we will give more examples of differential equations defining $\D$. Those examples will be used to characterise definable and algebraic closures of generic elements in the reducts. At the end of the section we will give two non-examples. Note that the results of this section will not be used later.

{We will show first that differential rational functions define the derivation.
\begin{proposition}\label{rational}
If $E(x,y)$ is given by $g(x) \cdot y=f(x)$ where $\frac{f(X)}{g(X)}$ is a differential rational function which is not an algebraic rational function, then $\D$ is definable in $\mathcal{F}_E$.
\end{proposition}}

{\begin{lemma}\label{two-derivations}
Let $\D_1$ and $\D_2$ be derivations on a field $K$ and $t \in K$ be such that $\D_1t=\D_2t=1$. If there is a non-zero algebraic polynomial $P(X_0,\ldots,X_n,Y_1,\ldots,Y_m)$ over $K$ such that  
$$P(X,\D_1 X,\ldots,\D_1^n X,\D_2 X,\ldots,\D_2^m X)=0$$
then $\D_1=\D_2$.
\end{lemma}
\begin{proof}
We can assume without loss of generality that $n=m$. As in the proof of Lemma \ref{L3-15} we can show there is a non-zero polynomial $P_1(\bar{X},\bar{Y})$ such that
$$P_1(X_1,\ldots,X_n, \D_2 Y - \D_1 Y + X_1, \ldots, \D^n_2 Y - \D^n_1 Y + X_n)=0.$$
Clearly $\D := \D_2-\D_1$ is a derivation of $K$. 
The above identity implies that for some non-zero polynomial $Q$ we have
$$Q(\D X, \D^2 X, \ldots, \D^n X)=0.$$
If $\D_1 \neq \D_2$ then $\D \neq 0$ and there is an element $b \in K$ with $\D b \neq 0$. Dividing $\D$ by $\D b$ we can assume that $\D b =1$. But then substituting $X\mapsto X+rb^j$ for $r \in \mathbb{Q}$ and $j=1, \ldots, n,$ we see that $Q=0$, which is a contradiction.
\end{proof}}

{\begin{proof}[Proof of Proposition \ref{rational}]
Suppose $$f(X)=P(X,\D X,\ldots,\D^nX),~ g(X)=Q(X,\D X, \ldots,\D^mX).$$ We will use Beth's definability theorem to show that $\D$ is definable in $T_E:=\Th(\mathcal{F}_E)$. Indeed, if we have two derivations $\D_1$ and $\D_2$ on a model $\mathcal{K}_E \models T_E$ that are compatible with $E$ (and $K$ is differentially closed with either of these derivations and $\D_1 t =\D_2t=1$), then
\begin{align*} \label{eq5-1}
& P(X,\D_1X,\ldots,\D_1^nX)\cdot Q(X,\D_2X,\ldots,\D_2^mX)=    \\
& P(X,\D_2X,\ldots,\D_2^nX)\cdot Q(X,\D_1X,\ldots,\D_1^mX).
\end{align*}
Since $f(X)/g(X)$ is not an algebraic rational function, the above identity shows that the conditions of Lemma \ref{two-derivations} are satisfied. Therefore $\D_1=\D_2$.
\end{proof}}

\begin{remark}
Note that even in the simple cases $y=\D^2 x$ and $y= (\D x)^2$ the differentiation is not definable without using $t$ since we can not distinguish between $\D$ and $-\D$.
\end{remark}

{Now we prove that if $E(x,y)$ defines an algebraic function of $x,\D x,\ldots,\D^n x$, i.e. $E$ is given by an equation $f(x,y)=0$ with $\ord_Y(f)=0$, then one can define $\D x$. But first we need to exclude some trivial counterexamples like $y \cdot \D x =0$ (see Example \ref{one-variable}).}

{\begin{definition}
A differential polynomial $f(X,Y)$ is said to be \textit{non-degenerate} if it cannot be decomposed into a product $g(X)h(X,Y)$ where $g$ is a differential polynomial and $h$ is an algebraic polynomial. An irreducible non-algebraic polynomial which depends on both variables is obviously non-degenerate.
\end{definition}}

{\begin{proposition}\label{algebraic}
Suppose $E(x,y)$ is defined by a non-degenerate equation $f(x,y)=0$ where $\ord_X(f)>0$ and $\ord_Y(f)=0$. Then $\D$ is definable in $\mathcal{F}_E$.
\end{proposition}}
\begin{proof}
{Pick a differentially transcendental element $a \in F$ and let $$f(a,Y)=\prod_{i=1}^k f_i(a,Y)^{e_i}$$ be the irreducible factorisation of $f(a,Y)$ over $k_0\langle a \rangle$. Denote $$g(a,Y):=\prod_{i=1}^k f_i(a,Y)=\sum_{i=0}^m g_i(a)\cdot Y^i,$$
where $g_i(X) \in k_0 \langle X\rangle$ and $g_m \neq 0$.} 

{Consider the formula
$$\psi(x,z_0,\ldots,z_m)= \exists y_1,\ldots,y_m \left(\bigwedge_{i \neq j} y_i \neq y_j \wedge \bigwedge_{i=1}^m E(x,y_i) \wedge \bigwedge_{i=1}^m \sum_{j=0}^m z_j\cdot y_i^j=0 \right).$$
Clearly, $E(a,y)$ holds if and only if $g(a,y)=0$. The polynomial $g(a,Y)$ has $m$ different roots. Therefore $\varphi(a,z_0,\ldots,z_m)$ holds if and only if the roots of $\sum_{i=0}^m z_i\cdot Y^i$ are exactly the same as those of $g(a,Y)$ (as these two polynomials have the same degree in $Y$). This can happen if and only if $\sum_{i=0}^m z_i\cdot Y^i$ is equal to $g(a,Y)$ up to a constant which depends on $a$. This means that
$$ \frac{z_i}{z_m} = \frac{g_i(a)}{g_m(a)}, $$
for all $i$. At least one of $\frac{g_i(X)}{g_m(X)}$ is not an algebraic rational function since otherwise $f$ would be degenerate. But then we can define $\D a$ from that differential rational function by Proposition \ref{rational} and we are done.}
\end{proof}

Next, we will apply Proposition \ref{algebraic} to work out definable and algebraic closures of generic points in the reducts. As before, let $a\in F$ be a generic point. We will show that the definable closure of $a$ in $\mathcal{F}_{R}$ coincides either with the definable closure in the differentially closed field or with that in the pure algebraically closed field.

It is well known what the definable and algebraic closures of arbitrary sets in differentially closed fields look like. Taking into account the fact that we have added $t$ as a constant symbol to the language, we see that for a set $A\subseteq F$ the definable and algebraic closures in $\mathcal{F}$ are given by $\dcl_{\D}(A)=k_0\langle A \rangle $ and $\acl_{\D}(A)=(k_0\langle A \rangle) ^{\alg} $, where $k_0=\mathbb{Q}(t)$ and $k_0\langle A \rangle$ is the differential subfield generated by $k_0$ and $A$. This immediately implies that in the reduct we have $k_0(A) \subseteq \dcl_{R}(A)\subseteq k_0\langle A \rangle $ and $(k_0(A))^{\alg} \subseteq \acl_{R}(A) \subseteq (k_0\langle A \rangle) ^{\alg} $.

We show that for generic elements one of these two extremal cases must happen.

\begin{theorem}\label{dcl,acl}
For $a\in F$ a generic point exactly one of the following statements holds:
\begin{itemize}
  \item $\dcl_{R}(a)=k_0(a)$; this holds if and only if $\acl_{R}(a)=(k_0(a))^{\alg}$ if and only if $\D$ is not definable;
  \item $\dcl_{R}(a)=k_0\langle a \rangle$; this holds if and only if $\acl_{R}(a)=(k_0\langle a \rangle)^{\alg}$ if and only if $\D$ is definable.
\end{itemize}
\end{theorem}
\begin{proof}
It will be enough to show that if $\acl_{R}(a)\supsetneq (k_0(a))^{\alg}$ then $\D$ is definable. Thus, let $\acl_{R}(a)\supsetneq (k_0(a))^{\alg}$. Choose $b\in (k_0\langle a \rangle)^{\alg} \setminus (k_0(a))^{\alg}$ which is algebraic (in the model theoretic sense) over $a$ in $\mathcal{F}_{R}$. There is a formula $\varphi(x,y) \in \mathfrak{L}_{R}$ such that $\varphi(a,b)$ holds and $\varphi(a,y)$ has finitely many realisations. Because $\varphi(a,y)$ defines a finite set in the differentially closed field $\mathcal{F}$, it is equivalent to an algebraic polynomial equation over $k_0 \langle a \rangle$. The latter is clearly non-degenerate and is not defined over $k_0(a)$ since $b$ is its root. Applying Proposition \ref{algebraic} we define $\D a$ (over $a$). Hence $\D$ is definable.
\end{proof}

Now using Proposition \ref{rational} we generalise Theorem \ref{finiterank}.

{\begin{proposition}\label{omega(n+1)}
Let $a \in F$ be a differentially transcendental element. If
$$\MR_R(a,\D a, \ldots, \D^n a) < \omega \cdot (n+1)$$
for some $n$ then $\D$ is definable.
\end{proposition}}
\begin{proof}
{We proceed to the proof by induction on $n$. The case $n=1$ is done in Theorem \ref{curve-def}. Assuming the theorem is true for all numbers less than $n$, we prove it for $n$.}

{There is a formula $\varphi(x_0,x_1,\ldots,x_n) \in \mathfrak{L}_R$ with $\MR_R(\varphi) < \omega \cdot (n+1)$ and
$$\mathcal{F}_R \models \varphi(a,\D a, \ldots, \D^n a).$$
Since $\varphi$ does not have ``full'' rank, we can assume without loss of generality it is given by a differential equation $f(x_0,x_1,\ldots,x_n)=0$ in the language of differential rings. Since $f(a,\D a, \ldots, \D^n a)=0$ and $a$ is generic, $f$ must be equal to $g(X,U_1,\ldots,U_n)$ for some differential polynomial $g$ with $g(X,0,\ldots,0)=0$ where $X:=X_0,~ U_i:=X_i - \D^i X$. Further, applying the method of additive translates as in the proof of Theorem \ref{curve-def}, we can assume that $g$ does not depend on the first variable, so we write $g(U_1,\ldots,U_n)$. However, we cannot proceed as in Theorem \ref{curve-def} and use multiplicative translates as there exist non-algebraic ``homogeneous'' differential polynomials of several variables.}

{\begin{claim}
$\D^n a \in \dcl_R(a,\D a, \ldots, \D^{n-1} a).$
\end{claim}}
\begin{proof}
{The set defined by $\varphi(a,\D a, \ldots, \D^{n-1} a, y)$ contains $\D^n a$. Moreover, that formula is given by $h(y-\D^n a)=0$ where $h(U)=g(0,\ldots,0, U)$. Consider the formula
\begin{equation}\label{big-formula}
\varphi\left[p(t)a, \D (p(t)a), \ldots, \D^{n-1}(p(t)a), \D^n(p(t)a)-p(t)\D^n a + p(t)y\right],
\end{equation}
for a non-zero polynomial $p(t) \in \mathbb{Q}[t]$.}

{It is easy to see that this is a formula in the language of reducts with parameters $a, \D a, \ldots, \D^{n-1}a$ and it is true of $y=\D^n a$. The formula \eqref{big-formula} is equivalent to $h(p(t)(y-\D^n a))=0$. Taking the conjunction of $\varphi(a,\D a, \ldots, \D^{n-1} a, y)$ and the formula \eqref{big-formula} we get a formula in the language $\mathfrak{L}_R$ equivalent to\footnote{Here we use square brackets for ease of reading. They do not have any special meaning.}
$$h[y-\D^n a]=0 \wedge h[p(t)(y-\D^n a)]=0.$$
This contains the point $\D^n a$ since $h(0)=0$ and is contained in sets defined by $h[p(t)(y-\D^n a)] - (p(t))^{\alpha}h[y-\D^n a]=0$ for $\alpha$ a positive integer. By an appropriate choice of $\alpha$ we reduce the multi-degree of $h$ and by a choice of $p$ we make sure the difference is not identically zero (see the proof of Theorem \ref{curve-def}). This can be done unless $h$ is an algebraic homogeneous polynomial. Iterating this process we will eventually reach a situation where $h$ has been replaced by an algebraic homogeneous polynomial in which case our formula defines $\D^n a$.}
\end{proof}

{We proved that there is a formula $\psi(\bar{x}) \in \mathfrak{L}_R$ such that $\psi(a,\D a, \ldots,\D^{n-1} a, y)$ has a unique solution which is $\D^n a$. We can assume $\MR_R(a,\D a, \ldots, \D^{n-1}a)= \omega \cdot n$ since otherwise $\D$ is definable by the induction hypothesis. Let $b_1, \ldots, b_{n-1}$ be differentially independent elements over $a$. Then
$$\tp_R(a,b_1,\ldots,b_{n-1}) = \tp_R(a,\D a, \ldots, \D^{n-1} a) =: p(x_0,\ldots,x_{n-1}).$$
Evidently $\exists ! y \psi(x_0, \ldots, x_{n-1},y) \in p$ where ``$\exists !$'' stands for ``there is a unique'' (it is obviously first-order expressible). Therefore
$$\mathcal{F}_R \models \exists ! y \psi(a,b_1 \ldots, b_{n-1},y),$$
and the unique solution of $\psi(a,b_1 \ldots, b_{n-1},y)$ is a differential rational function of $a,b_1,\ldots,b_{n-1}$. Denote it by $r(a,b_1,\ldots,b_{n-1})$. If $r$ is an algebraic rational function then
$$\xi(\bar{x}):= \forall y (\psi(x_0,\ldots,x_{n-1},y) \leftrightarrow y=r(x_0,\ldots,x_{n-1}))$$
is a formula in the language of reducts and is true of $(a,b_1,\ldots,b_{n-1})$. Hence it must be true of $(a,\D a, \ldots, \D^{n-1}a)$ too, which means $\D^n a = r(a,\D a, \ldots, \D^{n-1}a)$ which is impossible since $a$ is differentially transcendental.}

{Thus, $r$ is not algebraic. By a compactness argument (as in Lemma \ref{generic-withoutparameters}) we can choose $t_1,\ldots,t_{n-1} \in \mathbb{Q}(t,a)$ such that $\mathcal{F} \models \xi(a,t_1 \ldots, t_{n-1})$ and $r(a,t_1,\ldots,t_{n-1}) \in \mathbb{Q}(t)\langle a \rangle \setminus \mathbb{Q}(t,a)$. This guarantees that the formula $\psi(a,t_1 \ldots, t_{n-1},y)$ (which is in the language of reducts) defines a non-algebraic (in the field theoretic sense) element over $a$ and so $\dcl_R(a) \supsetneq k_0(a)$. So $\D$ is definable due to Theorem \ref{dcl,acl}.}
\end{proof}

Recall that in a stable theory a set $A$ (in the monster model) is called \textit{independent} (over $B$) if for any $a \in A$ we have $a \forkind_B A \setminus \{ a\}$.

{\begin{corollary}\label{independent}
$\D$ is definable in $\mathcal{F}_{R}$ if and only if the sequence $a, \D a, \D^2 a, \ldots $ is not independent (over the empty set) in $\mathcal{F}_{R}$.
\end{corollary}
\begin{proof}
If the sequence $a, \D a, \D^2 a, \ldots $ is not independent then for some $n$ the set $\{ a, \D a, \ldots, \D^n a\}$ is not independent. Therefore $\MR_R(a,\D a, \ldots, \D^n a) < \omega \cdot (n+1)$.
\end{proof}}

{As a common generalisation of Theorem \ref{curve-def} and Proposition \ref{algebraic} we prove the following result.}

\begin{proposition}
Suppose $E$ (a curve in general sense) contains a differential curve defined by a non-degenerate equation $f(x,y)=0$ where $\ord_X(f)>0$ and $\ord_Y(f)=0$. Then $\D$ is  definable in $\mathcal{F}_E$.
\end{proposition}
\begin{proof}
Let $g(X,Y)=p(X, \D X, \ldots, \D^n X, Y)$ be an irreducible non-degenerate factor of $f(X,Y)$. Furthermore, as $\ord_X(f)>0$ we can assume that $\ord_X(g)>0$. Consider the formula
$$\varphi(x,y_1,\ldots,y_n): = \exists z (E(x,z) \wedge p(x,y_1,\ldots,y_n,z)=0).$$
Clearly $\mathcal{F}_{E} \models \varphi (a, \D a, \ldots, \D^n a)$. Further, if $\varphi(a,b_1,\ldots, b_n)$ holds then for some $c$ we have $$p(a, \overline{b}, c) = 0 \wedge E(a,c).$$
Since $p$ is irreducible, $a, b_1,\ldots, b_n$ are algebraically dependent over $c$. Moreover, $\ord_X(g)>0$ implies that $b_1,\ldots, b_n$ are algebraically dependent over $\{a,c\}$. On the other hand, $c$ is differentially algebraic over $a$. Therefore $a, \bar{b}$ are differentially dependent and hence $\MR_{\D}(\varphi) < \omega \cdot (n+1)$. Now Proposition \ref{omega(n+1)} finishes the proof.
\end{proof}

One will certainly notice at this point that we found a number of conditions on $\mathcal{F}_{R}$ which are all equivalent to definability of $\D$. We sum up all these conditions in the following theorem.

\begin{theorem}\label{equivalent-conditions}
For a differentially transcendental element $a\in F$  the following are equivalent:
\begin{enumerate}
\item $\D$ is definable in the reduct $\mathcal{F}_{R}$ without parameters,

\item $\MR_{R}(\D a/a) < \omega$,

\item $\MR_{R}(\D a/a) = 0$,

\item $\tp_{R}(\D a/a)$ forks over the empty set,

\item The sequence $(\D^n a)_{ n \geq 0 }$ is not (forking) independent,

\item $\dcl_{R}(a) \supsetneq k_0(a)$,

\item $\acl_{R}(a)\supsetneq (k_0(a))^{\alg}$,

\item Every model of $T^+_{R}$ is the $\mathfrak{L}_R$-reduct (with canonical interpretation) of a differentially closed field,

\item Every automorphism of $\mathcal{F}_{R}$ fixes $\mathbb{D}$ setwise.
\end{enumerate}
\end{theorem} 
\begin{proof}
We need only show $9 \Rightarrow 1$. Take any automorphism $\sigma$ of $\mathcal{F}_{R}$ which fixes $a$. It fixes $\mathbb{D}$ setwise, hence $(\sigma(a),\sigma(\D a)) \in \mathbb{D}$. This means $\sigma(\D a) = \D (\sigma a) = \D a$. Thus any automorphism of $\mathcal{F}_{R}$ fixing $a$ fixes $\D a$. Since $\mathcal{F}_{R}$ is saturated, $\D a$ is definable over $a$. Therefore $\D$ is definable.
\end{proof}

We conclude this section by giving examples of differential equations that do not define $\D$. 

\begin{example}\label{one-variable}
We will show that unary relations cannot define $\D$. 

Let $R$ consist of unary relations, i.e. definable subsets of $F$ (by quantifier elimination of $\DCF_0$ we may assume $R$ consists of sets of solutions of one-variable equations). Then $\D$ is not definable in $\mathcal{F}_R$.

Consider the differential closure $K$ of $k_0$ inside $\mathcal{F}$, that is, 
$$K = \{ d \in F: \DR(d) < \omega \}.$$
This is by definition a differentially closed field. Take a generic element $a \in F$, i.e. an element outside $K$. Let $L \supseteq K$ be the differential closure of $K\langle a \rangle$ inside $\mathcal{F}$. Further, denote $a_i = \D^i a,~ i\geq 0$ and let $A$ be a transcendence basis of $L$ over $K$ containing these elements (not differential transcendence basis, which would consist only of $a$).

Define a new derivation $\D_1$ on $L$ as follows. Set $\D_1 = \D$ on $K \cup A\setminus \{a_0, a_1\}$ and $\D_1 a_0 = a_2,~ \D_1 a_1=a_0$. This can be uniquely extended to a derivation of $L$. The field automorphism $\sigma \in \Aut(L/K)$ which fixes $A\setminus \{a_0, a_1\}$ and swaps $a_0$ and $a_1$ is in fact an isomorphism of differential fields $\mathcal{L}=(L;+,\cdot, \D)$ and $\mathcal{L}_1=(L;+,\cdot, \D_1)$. Therefore the latter is differentially closed.

Thus we have a field $L$ equipped with two different derivations $\D$ and $\D_1$ and $L$ is a differentially closed field with respect to each of them. Further, $K \subset L$ consists of all differentially algebraic elements in $\mathcal{L}$. Since $\mathcal{L}$ and $\mathcal{L}_1$ are isomorphic over $K$, the differential closure of $k_0$ in $\mathcal{L}_1$ is equal to $K$ as well. Therefore the interpretations of relation symbols for one-variable differential equations in $\mathcal{L}$ and $\mathcal{L}_1$ are contained in $K$. But $\D$ and $\D_1$ agree on $K$ and therefore those interpretations agree in $\mathcal{L}$ and $\mathcal{L}_1$. This shows that $\D$ is not definable in the structure $\mathcal{F}_R$.

\end{example}

\begin{example}
Now we give a more interesting example.

\begin{proposition}\label{exp1}
The exponential differential equation $\D y = y \D x$ does not define $\D$.
\end{proposition}
We show first that for a differential equation $E$ if $\D$ is definable in $T_{E}$ then $E$ is uniquely determined by $T_{E}$.

\begin{lemma}
If $\D$ is definable in $T_{E}$ then for any differential equation $E'(x,y)$
$$T_{E} = T_{E'} \Rightarrow E = E'.$$
\end{lemma}
\begin{proof}
Let $E$ be given by the equation $f(x,y)=0$. Since $\D$ is definable, the formula $\forall x,y (E(x,y) \leftrightarrow f(x,y)=0)$ (more precisely, its translation into the language of the reducts) is in $T_{E}$. In other words, the fact that $E$ is defined by the equation $f(x,y)=0$ is captured by $T_{E}$. Therefore if $E'$ has the same theory as $E$ it must be defined by the same equation $f(x,y)=0$.
\end{proof}

\begin{proof}[Proof of Proposition \ref{exp1}]
An axiomatisation of the complete theory of the exponential differential equation is given in \cite{Kirby-semiab}. One can deduce from the axioms that the equation $\D y = 2y \D x$ is elementarily equivalent to the exponential equation. But clearly those two equations define different sets in differentially closed fields. Hence the previous lemma shows that $\D$ is not definable if $E$ is given by $\D y = y \D x$.
\end{proof}
We will give another proof of Proposition \ref{exp1} in Section \ref{Model completeness}. 
\end{example}

\section{Model completeness}\label{Model completeness}
\setcounter{equation}{0}

{In Section \ref{generic points} we showed that if a formula $\varphi(x,y)$ defines a small set which contains the point $(a,\D a)$ for a differentially transcendental element $a$ then $\D$ is definable. Moreover, if $\varphi$ is existential then $\D$ is existentially definable. Recall that smallness of a set can be verified as follows: if $b$ is a generic (differentially transcendental) element over $a$, that is, $(a,b)$ is a generic pair (differentially independent), then $\varphi(x,y)$ defines a small set if and only if $\neg \varphi(a,b)$. Thus, instead of working with formulas defining $\D$ we can work with formulas $\varphi(x,y)$ with $\varphi(a, \D a) \wedge \neg \varphi(a,b)$. 
\begin{definition}
A formula $\varphi(x,y) \in \mathfrak{L}_{R}$ is a $\D$-\textit{formula} if $\mathcal{F} \models \varphi(a, \D a) \wedge \neg \varphi(a,b)$, where $(a,b)$ is a differentially independent pair.
\end{definition}
Here we worked over the empty set. In particular, $a$ is differentially transcendental over the empty set and the definitions that we consider are again over the empty set, i.e. without parameters. However, it is clear that we could in fact work over any set $A \subseteq F$. In this case we should let $a$ be differentially transcendental over $A$. If $\varphi(x,y)$ is a formula over $A$ such that $\varphi(a, \D a) \wedge \neg \varphi(a,b)$ holds where $b$ is differentially transcendental over $Aa$ (in this case we will say $\varphi$ is a $\D$-formula over $A$), then certainly $\D$ is definable over $A$. Moreover, if $\varphi(x,y)$ is existential then $\D$ is existentially definable over $A$. In this section we use this fact to prove that under a natural assumption, if $\D$ is definable then it is existentially definable.} 



As above $a \in F$ is a differentially transcendental element and $k_0= \mathbb{Q}(t) = \dcl_R(\emptyset)$ (recall that $t$ is an element with $\D t=1$).

{\begin{theorem}\label{inductive}
If $T_{R}$ is inductive (i.e. $\forall \exists$-axiomatisable) and defines $\D$ then it defines $\D$ existentially {and, therefore,} $T_{R}$ is model complete.
\end{theorem}}

{This is similar to Lindstr\"om's theorem in spirit stating that an inductive theory, which is categorical in some infinite cardinality, is model complete (see \cite{TZ}). We can also consider another formulation of Theorem \ref{inductive}: if $T_R$ is not model complete and is inductive, then it does not define $\D$. In general, $T_R$ is not expected to be model complete, so in inductive reducts definability of a derivation is expected to be rare.}

We now establish an auxiliary result which will be used in the proof of Theorem \ref{inductive}.
\begin{lemma}\label{univ-exist-def}
Let $\varphi(x,\bar{u}) \in \mathfrak{L}_{R}$ be a quantifier-free formula and $p(X,Y,\bar{U}) \in k_0[X,Y,\bar{U}]$ be an algebraic polynomial which is monic in the $Y$ variable. Denote
$$\chi(x,y) := \forall \bar{u} (\varphi(x,\bar{u}) \rightarrow p(x,y,\bar{u})=0).$$
If $\exists \bar{u} \varphi(a,\bar{u})$ and $\chi(a, \D a)$ hold then $\D$ is existentially definable.
\end{lemma}
\begin{proof}
Let the tuple $(b_1,\ldots,b_m,e_1,\ldots,e_s)$ be of maximal differential transcendence degree $m$ over $a$ such that $\mathcal{F}_{R} \vDash \varphi(a,b_1,\ldots,b_m,e_1,\ldots,e_s)$ and assume that $b_1,\ldots,b_m$ are differentially independent over $a$.

Consider the formula
$$\psi(x,y,\bar{z}) = \exists v_1,\ldots,v_s ( \varphi(x,\bar{z},\bar{v}) \wedge p(x,y,\bar{z},\bar{v})=0 ).$$
Clearly $\psi(a,\D a,\bar{b})$ holds. Moreover, if $\psi(a,d,\bar{b})$ holds for some $d$, then for some $d_1,\ldots,d_s$ we have
$$\mathcal{F}_{R} \vDash \varphi(a,\bar{b},\bar{d}),$$
which implies that $d_1,\ldots,d_s$ have finite rank over $\{a,b_1,\ldots,b_m\}$. Since $p$ is monic as a polynomial of $Y$ and $p(a,d,\bar{b},\bar{d})=0$, we conclude that $d \in (k_0(a,\bar{b},\bar{d}))^{\alg}$ and hence $d$ is not generic over $\{a,b_1,\ldots,b_m\}$.

Thus working over the parameter set $B=\{b_1,\ldots,b_m\}$ we see that $a$ is generic over $B$ and $\psi(x,y,\bar{b})$ is a $\D$-formula over $B$. Hence we can make it into a proper definition of $\D$ with parameters from $B$. Thus, we get an existential definition of $\D$ with differentially independent parameters $b_1,\ldots,b_m$. By Lemma \ref{generic-withoutparameters} we have an existential definition without parameters.
\end{proof}

Now we are ready to prove our main theorem. 

\begin{proof}[Proof of Theorem \ref{inductive}]
Let $\delta(x,y)$ be a formula defining $\D$. We assume that $\D$ is not existentially definable, hence $\delta$ is not existential. The main idea of the proof is that unless one says explicitly that $\forall x \exists y \delta(x,y)$, one cannot guarantee that $\delta$ defines a function. In other words we will prove that $\forall x \exists y \delta(x,y)$ (which is not an $\forall \exists$-sentence) is not implied by the $\forall \exists$-part of $T_{R}$, as otherwise we will be able to find an existential definition of $\D$. This will contradict our assumption of inductiveness. 


Let $T$ be the $\forall \exists$-part of $T_R$, i.e. the subset of $T_R$ consisting of $\forall \exists$-sentences. In other words
$$T = \{ \forall \bar{x} \exists \bar{y} \varphi(\bar{x},\bar{y}): \varphi \mbox{ is a quantifier-free formula in } \mathfrak{L}_R,~ \mathcal{F}_R  \models \forall \bar{x} \exists \bar{y} \varphi(\bar{x},\bar{y}) \}.$$
Denote $\Phi := \{ \varphi(\bar{x},\bar{y}): \forall \bar{x} \exists \bar{y} \varphi(\bar{x},\bar{y}) \in T \}$.

By our assumption $T$ is an axiomatisation of $T_R$. However, we will get a contradiction to this by showing that $T$ has a model in which $\forall x \exists y \delta(x,y)$ does not hold. The construction of that model will go as follows. We start with the field $k=\mathbb{Q}(t,a)=k_0(a)$ and add solutions of the formulas $\varphi \in \Phi$ step by step (for $\varphi(\bar{x},\bar{y})\in \Phi$ we think of $\bar{x}$ as coefficients and of $\bar{y}$ as solutions). We also make sure that we do not add $\D a$ in any step. If the latter is not possible then we show that $\D$ is existentially definable.

In order to implement this idea, we expand the language by adding constant\footnote{In this proof we use the word ``constant'' for constant symbols only and not for constants in the sense of differential algebra. In particular, the interpretations of those constant symbols may not be constants in the differential sense.} symbols for solutions of all $\varphi \in \Phi$. First, take $C_0=\{a\}$. We will inductively add new constant symbols to $C_0$ countably many times.


If $C_l$ is constructed then $C_{l+1}$ is the expansion of $C_l$ by new constant symbols as follows. 
For each $\varphi(\bar{x}, \bar{y}) \in \Phi$ with $|\bar{x}|=m,~ |\bar{y}|=n$ say and for all $\bar{c} \in C_l^m$ add new constant symbols $d^1_{\varphi,\bar{c}},\ldots,d^n_{\varphi,\bar{c}}$. After adding these new constants for all $\varphi \in \Phi$ we get $C_{l+1}$. Finally set $C=\bigcup_lC_l$. This is a countable set. 

Now consider the following sets of sentences in the expanded language $\mathfrak{L}_R \cup C$. First, denote 
\begin{align*}
\Gamma(C) :=  \{ \varphi(c_1,\ldots,c_m,d^1_{\varphi,\bar{c}},\ldots,d^n_{\varphi,\bar{c}}) :
 \varphi(\bar{x},\bar{y}) \in \Phi,~ |\bar{x}|=m,~ |\bar{y}|=n,~  \bar{c} \in C^m\}.
\end{align*}
Furthermore, let
$$\Delta (C) := \{ \neg \delta(a,c): c \in C \}.$$
Finally we set 
$$\Sigma (C):= T_R \cup \tp_R(a) \cup \Gamma(C) \cup \Delta(C) .$$

\begin{claim}
$\Sigma:=\Sigma(C)$ is satisfiable.
\end{claim}
\begin{proof}
If it is not satisfiable, then a finite subset $\Sigma_0 \subseteq \Sigma$ is not satisfiable. Denote the set of constants from $C$ that occur in sentences from $\Sigma_0$ by $\{a,e_1,\ldots,e_n\}$ (if necessary, we can assume $a$ occurs in $\Sigma_0$ inessentially). We are going to give $a$ its canonical interpretation in $\mathcal{F}$ and this is the reason that we separated it from the other constant symbols. Let $\psi(a,e_1,\ldots,e_n) := \bigwedge (\Sigma_0 \cap \Gamma)$. The formula $\psi(x,u_1,\ldots,u_n)$ is clearly quantifier-free and without parameters.

Thus
$$ T_R \cup \tp_R(a) \cup\{ \psi(a,e_1,\ldots,e_n) \} \cup \{ \neg \delta(a,e_i): i=1,\ldots,n\}$$
is inconsistent. This means that in particular we cannot find interpretations for $e_1,\ldots,e_n$ in $\mathcal{F}_{R}$ which will make the latter into a model of $\Sigma_0$. As already mentioned above, $a$ is interpreted canonically in $\mathcal{F}$, i.e. its interpretation is the element $a \in F$.

Therefore
$$\mathcal{F}_{R} \nvDash \exists u_1,\ldots,u_n \left[ \psi(a,\bar{u}) \wedge \bigwedge_i \neg\delta(a,u_i) \right].$$

This means
$$\mathcal{F}_{R} \vDash \forall \bar{u} \left[ \psi(a,\bar{u}) \longrightarrow \bigvee_i u_i=\D a \right].$$

Note that evidently $\mathcal{F}_{R} \vDash \exists \bar{u} \psi(a,\bar{u})$, i.e. the implication above does not hold vacuously. So the formula 
$$\chi(x,y) := \forall \bar{u} \left[ \psi(x,\bar{u}) \longrightarrow \prod_{i=1}^n(y-u_i)=0 \right]$$
satisfies the conditions of Lemma \ref{univ-exist-def}. Hence, $\D$ is existentially definable. This contradiction proves the claim.
\end{proof}

Thus $\Sigma$ is satisfiable. Take a model $\mathcal{M}$ of $\Sigma$ and inside this model consider the subset $K$ consisting of interpretations of the constant symbols from $C$. We claim that $K$ is closed under addition and multiplication and contains $0,1,t$. This is because the sentence $\forall x, y \exists z,w (x+y=z \wedge x\cdot y = w)$, being $\forall \exists$, belongs to $T$. So, by our construction of $C$, for each $c_1, c_2 \in C$ we have elements $d_1, d_2 \in C$ such that the sentences $c_1+c_2=d_1,~ c_1 \cdot c_2 = d_2$ are in $\Sigma$. Similarly $0,1,t \in K$ since the sentences $\exists x (x=0),~ \exists x (x=1)$, and $\exists x(x=t)$ are in $T$. Therefore $K$ is a structure in the language of rings. In fact it is an algebraically closed field (containing $k$) since $\ACF_0$ is $\forall \exists$-axiomatisable. Hence $K$ is a structure $\mathcal{K}_{R}=(K;+,\cdot,0,1,t,P)_{P \in R}$ in the language of the reducts (with the induced structure from $\mathcal{M}$). By the choice of $\Sigma$ we know that $\mathcal{K}_{R}$ is a model of $T$. 

If we choose $\mathcal{M}$ to be saturated of cardinality $|F|$ (such a model exists due to stability) then we can identify it with $\mathcal{F}_R$.\footnote{Note that this is not essential, but it helps to understand how the proof works.} In that case $\mathcal{K}_R$ is a substructure obtained by starting with $k_0(a)$ and inductively adding solutions to formulas from $\Phi$. 

Suppose for a moment that $\delta$ is universal in order to illustrate what we are going to do next. Let $$\delta(x,y) = \forall \bar{v} \rho (x,y, \bar{v})$$ with $\rho$ quantifier-free. 
Since $\mathcal{F}_R \models \neg \delta (a,s)$ for any $s \in K$, there is a witness $\bar{l}_s \in F$ such that $\mathcal{F}_R \models \neg \rho (a,s, \bar{l}_s)$. However this witness may not be in $K$. So we add all those witnesses to $K$ and then repeat the above procedure to make it a model of $T$. We also make sure we never add $\D a$, which is possible as above (otherwise $\D$ would be existentially definable). Iterating this process countably many times and taking the union of all the constructed substructures we end up with a structure $\mathcal{N}_R$ in the language of reducts which is a model of $T$ and contains witnesses for each of the formulas $\exists  \bar{v}\neg \rho(a,s,\bar{v})$ where $s \in N$. Thus,  $\mathcal{N}_R \models \neg \exists y \delta(a,y)$ which means that $T$ is not an axiomatisation of $T_{R}$. This contradiction proves the theorem.

{Now we consider the general case. Let $\delta$ be of the form
$$\delta(x,y)=\forall \bar{v}_1 \exists \bar{w}_1 \forall \bar{v}_2 \ldots \forall \bar{v}_n \exists \bar{w}_n \rho (x,y,\bar{v}_1,\ldots,\bar{v}_n),$$
where $\rho$ is quantifier-free (the tuples $\bar{v}_1$ and $\bar{w}_n$ can be empty). Then
$$\neg \delta(x,y)=\exists \bar{v}_1 \forall \bar{w}_1 \exists \bar{v}_2 \ldots \exists \bar{v}_n \forall \bar{w}_n \neg \rho (x,y,\bar{v}_1,\ldots,\bar{v}_n).$$
We add new constant symbols as follows. Firstly, for each $s \in C$ we add a tuple of constants $\bar{l}^1_s$ of the same length as $\bar{v}_1$. Then for each $i$ and each tuple $\bar{c} \in C^{| \bar{w}_{i}|}$ we add new constants  $\bar{l}^{i+1}_{\bar{c}}$ with $|\bar{l}^{i+1}_{\bar{c}}|=|\bar{v}_{i+1}|$. Denote this extension of $C$ by $C'$. Then we add new constant symbols to $C'$ for solutions of all formulas $\varphi \in \Phi$ as above. We denote this set by $C^1$. Then we iterate this procedure by adding new constants to witness $\neg \delta(a,s)$ (for each $s$ from the set of constants already constructed) and then adding new constants for solutions of $\varphi \in \Phi$. Thus, we get a chain $C\subseteq C^1 \subseteq C^2 \subseteq \ldots$. Let $\tilde{C}$ be their union.\footnote{{In this construction for an element $s \in C$ (and for tuples $\bar{c}_i$) we add one $\bar{l}^1_s$ in each step. Though this does not cause any problems, at each step we can add the corresponding sets of constants only for new constant symbols. In particular, after adding one $\bar{l}^1_s$ we do not add any such tuple for the same $s$ any more. Alternatively, we could just require all those different tuples for the same element $s$ to be equal by adding the appropriate formulas stating their equality.}}}

{For $A \subseteq \tilde{C}$ denote $$\Xi(A):=\{ \neg \rho(a,s,\bar{l}^1_s,\bar{c}_1,\bar{l}^2_{\bar{c}_1},\ldots,\bar{l}^n_{\bar{c}_{n-1}},\bar{c}_n) : s \in A,~ \bar{c}_i \in A^{| \bar{w}_i|} \}.$$
For any $s \in K$ we know that $\mathcal{F}_R \models \neg \delta (a,s)$, therefore
$\Sigma(C)\cup \Xi(C)$ is satisfiable (note that this collection of sentences contains parameters from $C'$). The proof of the above claim shows that $\Sigma(C^1) \cup \Xi(C)$ is satisfiable and so $\Sigma(C^1) \cup \Xi(C^1)$ is satisfiable too. Proceeding inductively we see that $\Sigma(C^i) \cup \Xi(C^i)$ is satisfiable for each $i < \omega$. Hence, by compactness, $\Sigma(\tilde{C}) \cup \Xi(\tilde{C})$ is satisfiable.}

{The interpretation of $\tilde{C}$ in a model of $\Sigma(\tilde{C}) \cup \Xi(\tilde{C})$ gives a structure $\mathcal{N}_R$ in the language of reducts which is a model of $T$ and contains witnesses for each of the formulas $\exists \bar{v}_1 \forall \bar{w}_1 \exists \bar{v}_2 \ldots \exists \bar{v}_n \forall \bar{w}_n \neg \rho (a,s,\bar{v}_1,\ldots,\bar{v}_n)$ where $s \in N$. Hence  $\mathcal{N}_R \models \neg \exists y \delta(a,y)$ which means that $T$ is not an axiomatisation of $T_{R}$, which is a contradiction.}
\end{proof}

As an immediate application of Theorem \ref{inductive} we give another proof of Proposition \ref{exp1} which states that if $E$ is the exponential differential equation, i.e. it is given by $\D y = y \D x$, then $\D$ is not definable in $\mathcal{F}_{E}$.\footnote{Note that Kirby \cite{Kirby-semiab} gives yet another proof of this fact by considering a lattice of reducts for exponential differential equations of some collections of semiabelian varieties, and showing that each of these can be expanded properly inside $\DCF_0$.} Indeed, Kirby gives an $\forall \exists$-axiomatisation of the first-order theory of the exponential differential equation (\cite{Kirby-semiab}). It is not model complete however, hence $\D$ cannot be definable due to Theorem \ref{inductive}. Of course it is the Ax-Schanuel inequality that is responsible for this. As Kirby proved, it is an adequate predimension inequality.

\addcontentsline {toc} {section} {Bibliography}
\bibliographystyle {alpha}
\bibliography {ref}

\end{document}